\documentclass{article} 

\usepackage{nips10submit_e,times}
\usepackage{amsmath}
\usepackage[algoruled,noend]{algorithm2e}
\usepackage{hyperref}
\usepackage[dvips]{graphicx}
\usepackage[utf8]{inputenc}
\usepackage{subfigure}

\title{Fast and Faster: A Comparison of Two Streamed Matrix Decomposition Algorithms}

\author{
Radim Řehůřek\\
NLP lab, Masaryk University in Brno, Czech Republic\\
\texttt{radimrehurek@seznam.cz}
}

\nipsfinalcopy 

\begin{document}

\maketitle

\begin{abstract}
With the explosion of the size of digital dataset, the limiting factor for decomposition algorithms is the \emph{number of passes} over the input, as the input is often stored out-of-core or even off-site. Moreover, we're only interested in algorithms that operate in \emph{constant memory} w.r.t. to the input size, so that arbitrarily large input can be processed. In this paper, we present a practical comparison of two such algorithms: a distributed method that operates in a single pass over the input vs. a streamed two-pass stochastic algorithm. The experiments track the effect of distributed computing, oversampling and memory trade-offs on the accuracy and performance of the two algorithms. To ensure meaningful results, we choose the input to be a real dataset, namely the whole of the English Wikipedia, in the application settings of Latent Semantic Analysis.
\end{abstract}

\section{Introduction}

Matrix decomposition algorithms are commonly used in a variety of domains across much of the field of Computer Science\footnote{Examples include Latent Semantic Analysis in Natural Language Processing; (discrete) Karhunen–Loève Transform in Image Processing or Recommendation Systems in Information Retrieval. SVD is also used in solving shift-invariant Differential Equations, in Geophysics, in Signal Processing, in Antenna Array Processing, \ldots}. 
Research has traditionally focused on optimizing the number of FLOPS (floating point operations) and numerical robustness of these algorithms~\cite{comon1990tracking,golub1996matrix}. However, modern datasets are too vast to be stored in main memory, or even on a single computer, so that communication itself quickly becomes a bottleneck.

One of the oldest and most widely known matrix decomposition algorithms is the Singular Value Decomposition (SVD), or its closely related eigen decomposition, which produce a provably optimal (in the least-squares sense) rank-$k$ factorizations when truncated. In the following, $n$ will denote the number of observations (matrix columns), $m$ the number of features (matrix rows) and $k$ the truncated target rank, $k \ll m \ll n$. In practise, the optimal decompositions are notoriously expensive to compute and truly large-scale applications are rare. The most common remedy is a) approximation (subsampling the input), b) some sort of incremental updating scheme which avoids recomputing the truncated models from scratch every time an observation/feature is updated, or c) giving up on a globally optimal solution and using another, heuristic algorithm. 
One way or another, the algorithm must avoid asking for $O(n)$ memory, as the number of observations is assumed to be too large in modern problems. Table~\ref{tab:algos} summarizes available algorithms (and their implementations) with respect to several interesting characteristics, such as whether or not they are distributed, whether they can be incrementally updated, how many input passes are required or whether they realize subspace tracking (infinite input stream, gradual model decay).

This paper compares two particular modern approaches to large-scale eigen decomposition: a one-pass streamed distributed algorithm from~\cite{lsa_ecir} and a modified stochastic streamed two-pass algorithm from~\cite{halko2009finding}. They require one and two passes over the input respectively; we will call them $P1$ and $P2$ from now on. Both are streamed, meaning no random access to observations is required and their memory requirements are constant in the number of observations. Some modifications to the original $P2$ algorithm were necessary to achieve this; these are described below. Apart from the practical side-by-side comparison, we also present a hybrid of the two methods here, a novel algorithm which takes advantage of the speed of $P2$ while retaining the one-pass quality of $P1$. 

\begin{table*}[tb]
\renewcommand{\baselinestretch}{1.2}

\centering
\scriptsize
\caption{Selected algorithms for truncated, partial eigen-decomposition and their characteristics. ``---'' stands for \emph{no/not found}.}
\label{tab:algos} 
\centering

\begin{tabular}{ p{0.20\textwidth} |  c  | c|  c | c | c | p{0.25\textwidth}}
  \hline
  \textbf{Algorithm} & \multicolumn{1}{c}{\textbf{Distributed}} & \multicolumn{2}{c}{\textbf{Incremental in}} & \textbf{\# passes} & \textbf{Subspace} & \textbf{Implementations} \\
  &  & observations & features &  & \textbf{tracking} &  \\
  \hline
	Krylov subspace methods (Lanczos, Arnoldi) & yes & --- & --- & $O(k)$ & --- & \href{http://soi.stanford.edu/~rmunk/PROPACK/}{PROPACK}, \href{http://www.caam.rice.edu/software/ARPACK/}{ARPACK}, \href{http://www.netlib.org/svdpack/}{SVDPACK}, \href{http://mahout.apache.org/}{MAHOUT}, \ldots \\
 \hline
	\cite{halko2009finding} & yes & --- & --- & $O(1)$ & --- & \href{http://code.google.com/p/redsvd/}{redsvd}, \href{http://cims.nyu.edu/~tygert/pca.m}{pca.m}, our own \\
 \hline
	\cite{gorrell2005generalized} & --- & --- & --- & $O(k)$ & --- & \href{http://alias-i.com/lingpipe/}{LingPipe}, our own \\
 \hline
	\cite{zha1999updating} & --- & yes & yes & 1 & yes & ---, our own \\
\hline
	\cite{levy2000sequential} & --- & yes & --- & 1 & yes & ---, our own \\
 \hline
 	\cite{brand2006fast} & --- & yes & yes & 1 & --- & ---, our own \\
 \hline
  	\cite{lsa_ecir} & yes & yes & --- & 1 & yes & \href{http://nlp.fi.muni.cz/projekty/gensim/}{our own, open-sourced}\\
  \hline
\end{tabular}
\vspace{-0.3cm}
\renewcommand{\baselinestretch}{1.0}
\end{table*}

\subsection{Stochastic two-pass algorithm, $P2$}

The one-pass stochastic algorithm as described in~\cite{halko2009finding} is unsuitable for large-scale decompositions, because the computation requires $O(nk+mk)$ memory. We can reduce this to a managable $O(mk)$, i.e. independent of the input stream size $n$, at the cost of running two passes over the input matrix instead of one\footnote{Actually, $2+q$ passes are needed when using $q$ power iterations.}. This is achieved by two optimizations: 1) the sample matrix is constructed piece-by-piece from the stream, instead of a direct matrix multiplication, and 2) the final dense decomposition is performed on a smaller $k \times k$ eigenproblem $BB^T$ instead of the full $k \times n$ matrix $B$. 

These two ``tricks'' allow us to compute the decomposition in constant memory, by processing the observations one after another, or, preferrably, in as large chunks as fit into core memory.
The intuition behind these optimizations if fairly straightforward, so we defer fleshing out the full algorithm to Appendix~\ref{app:p2}.

\subsection{One-pass algorithm, $P1$}

Streamed one-pass algorithms are fundamentally different from the 2-pass algorithm above (or any other multi-pass algorithm), in that as long as they manage to keep their memory requirements constant, they allow us to process infinite input streams. In environments where the input cannot be persistently stored, this may be the only option. 

In~\cite{lsa_ecir}, I describe one such algorithm. It works by computing in-core decompositions of document chunks, possibly on different machines, and efficiently merging these dense partial decompositions into one. The partial in-core decomposition algorithm is viewed as ``black box'' and chosen to be Douglas Rohde's \href{http://tedlab.mit.edu/~dr/SVDLIBC}{SVDLIBC}. The coarsely-grained parallelism of this algorithm makes it suitable for distributing the computation over a cluster of commodity computers connected by a high-latency network.

\subsection{Hybrid algorithm, $P12$}

In this work, we also explore combining the two above approaches. We consider using the in-core stochastic decomposition of~\cite{halko2009finding} instead of SVDLIBC in the one-pass merging framework of~\cite{lsa_ecir}. This hybrid approach is labelled $P12$ in the experiments below.

\section{Experiments}

We will be comparing the algorithms on an implicit 100,000$\times$ 3,199,665 sparse matrix with 0.5 billion non-zero entries (0.15\% density). This matrix represents the entire English Wikipedia\footnote{Static dump as downloaded from~\url{http://download.wikimedia.org/enwiki/latest}, June 2010.}, with the vocabulary (number of features) clipped to the 100,000 most frequent word types\footnote{The corpus preprocessing setup is described in more detail \href{http://nlp.fi.muni.cz/projekty/gensim/wiki.html}{online.}}. 
In all experiments, the number of requested eigen factors is arbitrarily set to $k=400$.

The experiments used three 2.0GHz Intel Xeon workstations with 4GB of RAM, connected by Ethernet on a single network segment. The machines were not dedicated; due to the large amount of experiments, we only managed to run each experiment twice. We report the better of the two times.

\subsection{Oversampling}
\label{exp1:text}

In this set of experiments, we examine the relative accuracy of the three algorithms. $P2$ has two parameters which affect accuracy: the oversampling factor $l$ and the number of power iterations $q$. In the one-pass algorithms $P1$ and $P12$, we improve accuracy by asking for extra factors $l$ during intermediate computations, to be truncated at the very end of the decomposition. 

\begin{figure*}[H]
\caption{Effects of the oversampling parameter $l$ on accuracy (Experiment~\ref{exp1:text}). Wall-clock times are in brackets. Experiments were run on a single machine, with chunks of 20,000 documents.}
\label{exp1:results}

\subfigure[Oversampling for $P1$, $P2$ and $P12$ algorithms.] {
	\includegraphics[scale=0.36,clip]{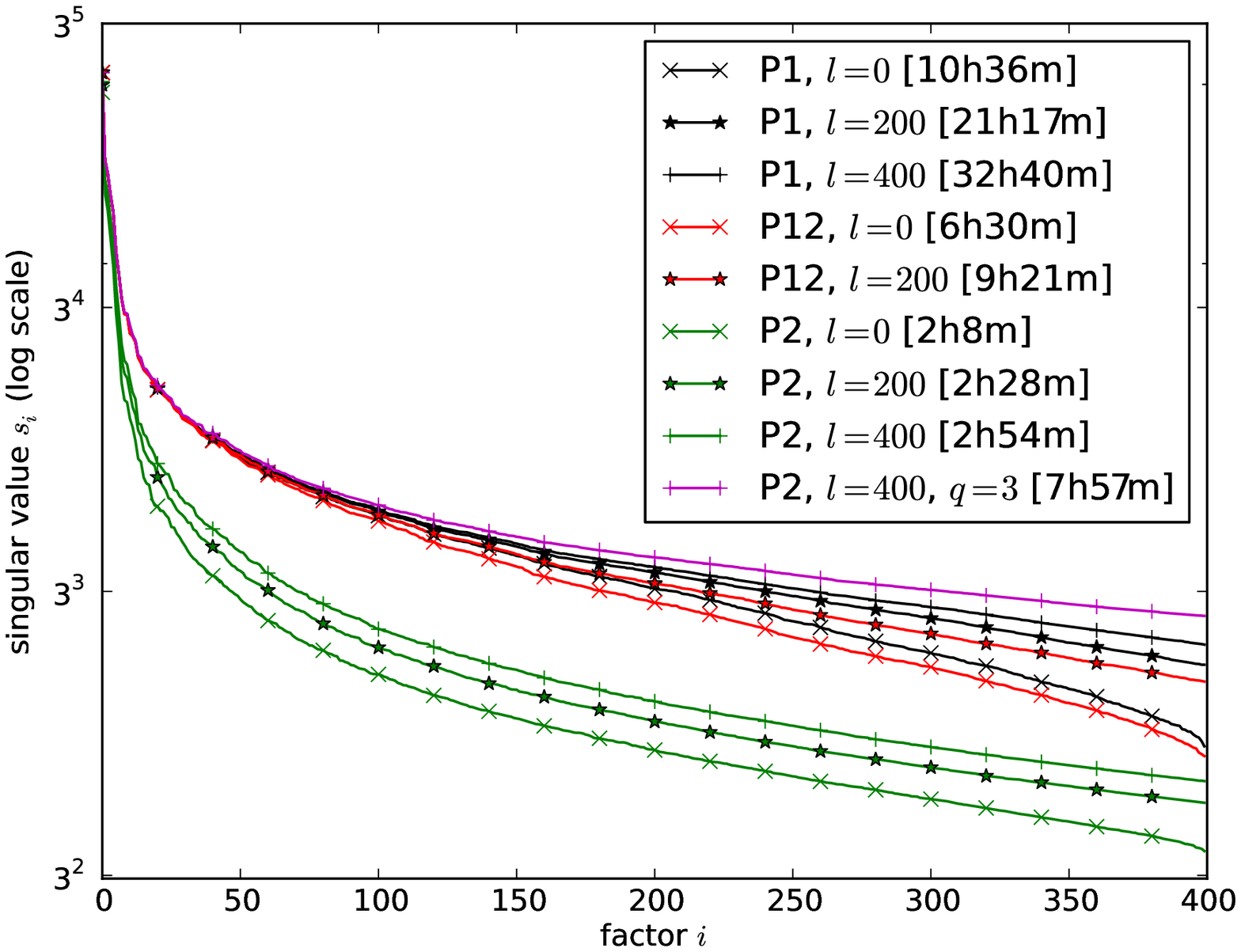}
}
\subfigure[Oversampling and power iterations for the $P2$ algorithm.] {
	\includegraphics[scale=0.36,clip]{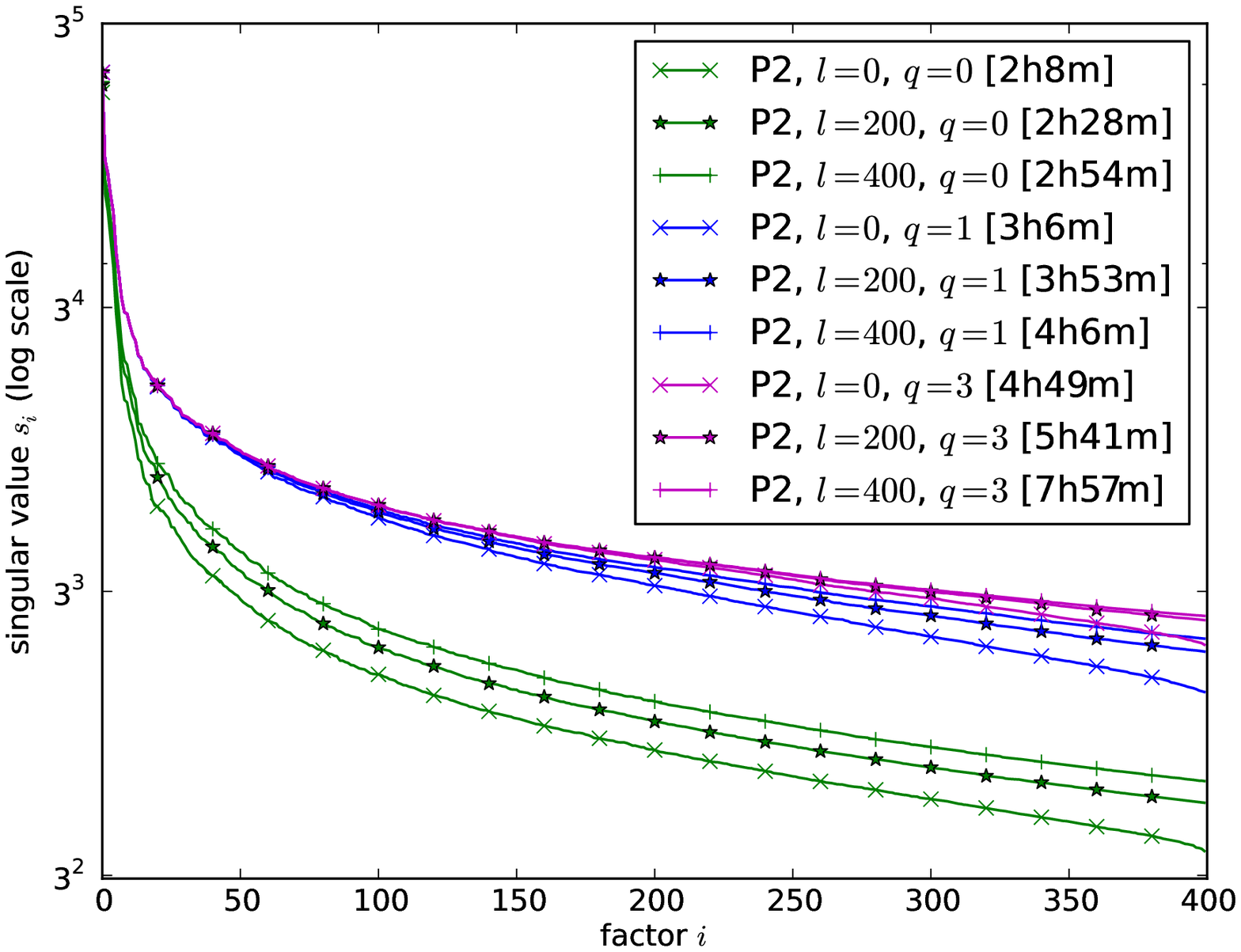}
}
\end{figure*}

Figure~\ref{exp1:results} summarizes both the relative accuracy and runtime performance of the algorithms, for multiple choices of $l$ and $q$. We see that although all methods are very accurate for the greatest factors, without oversampling the accuracy quickly degrades. This is especially true of the $P2$ algorithm, where no amount of oversampling helps and power iterations are definitely required. 

The ``ground-truth'' decomposition is unknown, so we cannot give absolute errors. However, according to our preliminary experiments on a smaller corpus, the stochastic algorithm with extra power iterations and oversampling gives the most accurate results; we will therefore plot it in all subsequent figures, in magenta colour, as a frame of reference. Note that all algorithm consistently err on the side of \emph{underestimating} the magnitude of the singular values---as a rule of thumb, the greater the singular values in each plot, the more accurate the result.

\subsection{Chunk size}
\label{exp2:text}

\begin{figure*}[H]
\begin{minipage}[b]{0.48\textwidth}
\centering
\caption{Accuracy and wall-clock times for different chunk sizes in $P1$ and $P12$ (Experiment~\ref{exp2:text}), no oversampling.}
\label{exp2:results}
\centering
\includegraphics[scale=0.36,clip]{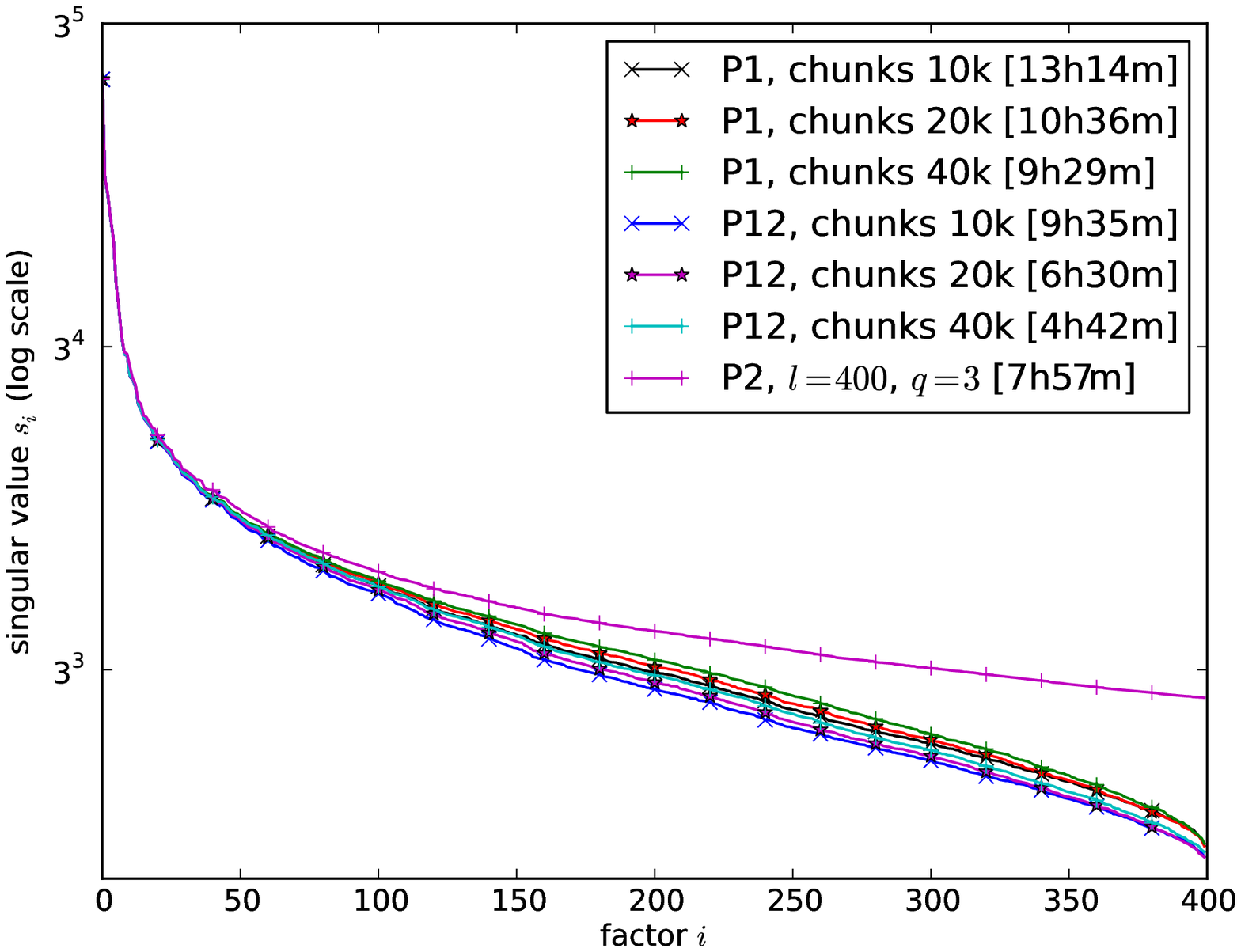}
\end{minipage}
\qquad
\begin{minipage}[b]{0.48\textwidth}
\centering
\caption{Effects of input order on the $P1$ algorithm (Experiment~\ref{exp3:text}). Chunk size is set to 40,000 documents, no oversampling.}
\label{exp3:results}

\centering
\includegraphics[scale=0.36,clip]{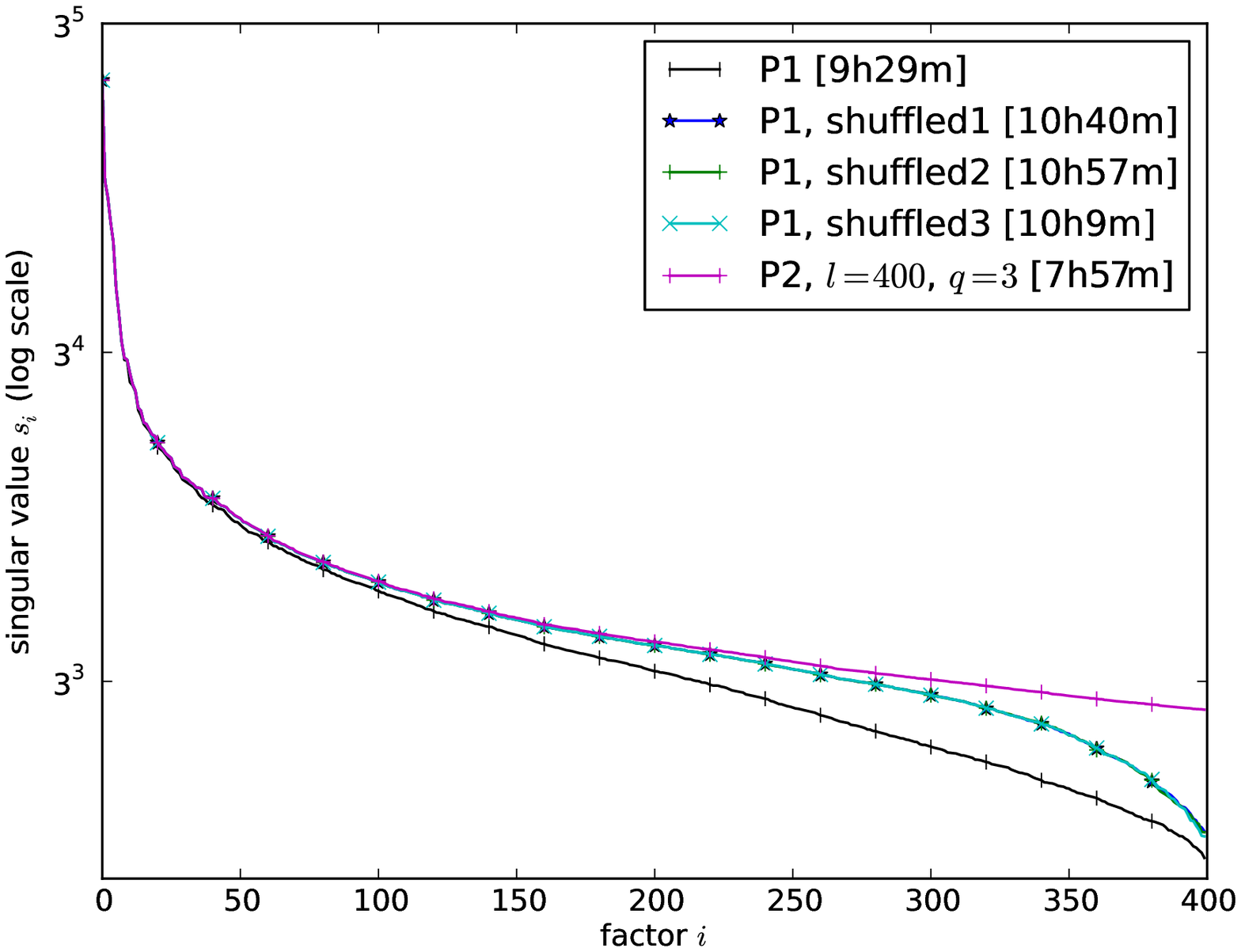}
\end{minipage}
\end{figure*}

The one-pass algorithms $P1$ and $P12$ proceed in document chunks that fit into core memory. A natural question is, what effect does the size of these chunks have on performance and accuracy? With smaller chunks, the algorithm requires less memory; with larger chunks, it performs fewer merges, so we might expect better performance. This intuition is quantified in Figure~\ref{exp2:results}, which lists accuracy and performance results for chunk sizes of 10,000, 20,000 and 40,000 documents. 

We see that chunk sizes in this range have little impact on accuracy, and that performance gradually improves with increasing chunk size. This speed-up is inversely proportional to the efficiency of the decomposition merge algorithm: with a hypothetical zero-cost merge algorithm, there would be no improvement at all, and runtime would be strictly dominated by costs of the in-core decompositions. On the other hand, a very costly merge routine would imply a linear relationship.

\subsection{Input stream order}
\label{exp3:text}

In the Wikipedia input stream, observations are presented in lexicographic order---observation corresponding to the Wikipedia entry on \emph{anarchy} comes before the entry on \emph{bible}, which comes before \emph{censorship} etc. This order is of course far from random, so we are naturally interested in how it affects the resulting decomposition of the single-pass algorithms (the two-pass algorithm is order-agnostic by construction).

To test this, we randomly shuffled the input stream and re-ran the experiments on $P1$. Ideally, the results should be identical, no matter how we permute the input stream. Results in Figure~\ref{exp3:results} reveal that this is not the case: singular values coming from the shuffled runs are distinctly different to the ones coming from the original, alphabetically ordered sequence. This likely shows that the one-pass truncated scheme has some difficulties adjusting to gradual subspace drift. With the shuffled input, no significant drift can occur thanks to the completely random observation order, and a much higher accuracy is retained even without oversampling.

\subsection{Distributed computing}
\label{exp4:text}

\begin{figure*}[H]
\begin{minipage}[b]{0.48\textwidth}
\caption{Distributed computing for algorithms $P1$, $P12$ (Experiment~\ref{exp4:text}). The chunk size is set to 20,000 documents, no oversampling.}
\label{exp4:results}
\centering
\includegraphics[scale=0.36,clip]{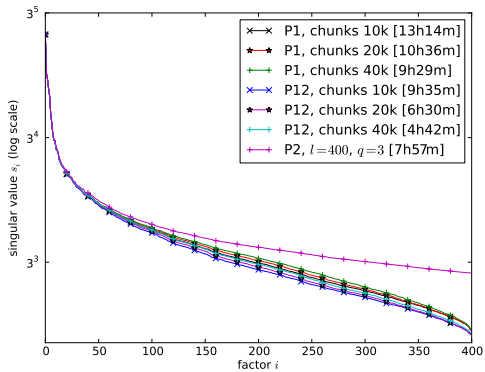}
\end{minipage}
\qquad
\end{figure*}

The two single pass algorithms, $P1$ and $P12$, lend themselves to easy parallelization. In Figure~\ref{exp4:results}, we evaluate them on a cluster of 1, 2 and 4 computing nodes. The scaling behaviour is linear in the number of machines, as there is virtually no communication going on except for dispatching the input data and collecting the results. As with chunk size, the choice of cluster size does not affect accuracy much.

The $P2$ algorithm can be distributed too, but is already dominated by the cost of accessing data in its $q+2$ passes. Routing data around the network gives no performance boost, so we omit the results from the figure. We note that distributing $P2$ would still make sense under the condition that the data is already predistributed to the computing nodes, perhaps by means of a distributed filesystem.


\section{Conclusion}

We presented a streamed version of a two-pass stochastic eigen decomposition algorithm and compared it to two streamed one-pass algorithms, one of which is a novel one-pass distributed algorithm. The comparison was done in the context of Latent Semantic Analysis, on a corpus of 3.2 million documents comprising the English Wikipedia.

On a single 2GHz machine, the top achieved decomposition times were 4 hours and 42 minutes for the one-pass $P12$ algorithm and 3 hours 6 minutes for the stochastic multi-pass algorithm. Without power iterations and with reduced amount of oversampling, we recorded even lower times, but at the cost of a serious loss of accuracy. On a cluster of four computing nodes on three physical machines, the single pass $P12$ decomposition was completed in 1 hour and 41 minutes.

We observed that the lightning-fast stochastic algorithm suffers from serious accuracy issues, which can be remedied by increasing the number of passes over the input (power iterations), as suggested in~\cite{halko2009finding}. But, as the number of passes is the most precious resource in streaming environments, the otherwise slower one-pass algorithms become quickly competitive.  The one-pass algorithms, one the other hand, suffer from dependency on the order of observations in the input stream; we will return to this behaviour in future work.

A practical and perhaps even more exciting contribution is a modern implementation of these algorithms that we release into open-source as \href{http://nlp.fi.muni.cz/projekty/gensim/}{\emph{gensim}}. Written in Python, it still manages to get top performance thanks to the use of Python's NumPy library with fast BLAS calls under the hood.

\subsubsection*{Acknowledgments}

This study was partially supported by the LC536 grant of MŠMT ČR.

\bibliographystyle{apalike}
{\scriptsize
\bibliography{rehurek_nips}}

\renewcommand{\baselinestretch}{1.0}

\clearpage

\appendix

\section{Streamed Stochastic Eigen Decomposition}
\label{app:p2}

\renewcommand{\baselinestretch}{1.8}

\smallskip
\begin{algorithm}[H]
\SetCommentSty{}

\BlankLine

\KwIn{$m \times n$ input matrix $A$, presented as a stream of observation chunks $A=\begin{bmatrix}C_1, C_2, \ldots, C_C\end{bmatrix}$. Truncation factor $k$. Oversampling factor $l$. Number of power iterations $q$.}
\KwOut{$U$, $S^2$ spectral decomposition of $A$ (i.e, $US^2U^T=AA^T$) truncated to the $k$ greatest factors.}
\KwData{Intermediate matrices require $O(m(k+l))$ memory; in particular, the algorithm avoids materializing any $O(n)$ or $O(m^2)$ matrices.}

\BlankLine
\BlankLine

\tcp{Construct the $m \times (k+l)$ sample matrix $Y=AO$, in one pass over the input stream.}
$Y \leftarrow$ \textbf{sum}($C_i O_i$ \textbf{for} $C_i$ \textbf{in} $A$) \tcp*[l]{each $O_i$ is a random $|C_i| \times (k+l)$ gaussian matrix}

\BlankLine

\tcp{Run $q$ power iterations to improve accuracy (optional), $Y=(AA^T)^qAO$. Needs $q$ extra passes.}
\For{$iteration \leftarrow 1$ \KwTo $q$}{
	$Y \leftarrow $ \textbf{sum}($C_i (C_i^T Y)$ \textbf{for} $C_i$ \textbf{in} $A$)\;
}

\BlankLine
\tcp{Construct the $m \times (k+l)$ orthonormal action matrix $Q$, in-core.}
$Q \leftarrow orth(Y)$\;

\BlankLine
\tcp{Construct $(k+l)\times(k+l)$ covariance matrix $X=BB^T$ in one pass, where $B=Q^TA$.}

$X \leftarrow $ \textbf{sum}($(Q^TC_i)(Q^TC_i)^T$ \textbf{for} $C_i$ \textbf{in} $A$) \tcp*[l]{BLAS rank-$k$ update routine SYRK}


\BlankLine
\tcp{Compute $U$, $S$ by means of the small $(k+l)\times(k+l)$ matrix $X$.}

$U_X, S_X \leftarrow eigh(X)$\;

\BlankLine
\tcp{Go back from the eigen values of $X$ to the eigen values of $B$ (= eigen values of $A$).}

$S^2 \leftarrow$ first $k$ values of $\sqrt{S_X}$\;
$U \leftarrow $ first $k$ columns of $Q U_X$\;

\caption{Two-pass Stochastic Decomposition in Constant Memory with Streamed Input}
\label{app:p2}
\end{algorithm}
\medskip

\clearpage

\section{Wikipedia LSA Topics}
\label{wiki_topics}

First ten topics coming from the $P2$ decomposition with three power iterations and 400 extra samples. The top ten topics are apparently dominated by meta-topics of Wikipedia administration and by robots importing large databases of countries, films, sports, music etc.

\renewcommand{\baselinestretch}{1.3}
\small

\bigskip

\begin{tabular}{ c | c | p{0.75\textwidth} }
  \hline
  \textbf{Topic $i$} & \multicolumn{1}{c}{\textbf{Singular}} & \textbf{Ten most salient words for topic $i$, with their weights} \\
  & \textbf{ value $s_i$} & \\
 \hline
 	1. & 201.118 & -0.474*``delete" + -0.383*``deletion" + -0.275*``debate" + -0.223*``comments" + -0.220*``edits" + -0.213*``modify" + -0.208*``appropriate" + -0.194*``subsequent" + -0.155*``wp" + -0.117*``notability" \\
 \hline
 	2. & 143.479 & 0.340*``diff" + 0.325*``link" + 0.190*``image" + 0.179*``www" + 0.169*``user" + 0.157*``undo" + 0.154*``contribs" + -0.145*``delete" + 0.116*``album" + -0.111*``deletion" \\
 \hline
 	3. &136.235 & 0.421*``diff" + 0.386*``link" + 0.195*``undo" + 0.182*``user" + -0.176*``image" + 0.174*``www" + 0.170*``contribs" + -0.111*``album" + 0.105*``added" + -0.101*``copyright" \\
 \hline
 	4. & 125.436 & 0.346*``image" + -0.246*``age" + -0.223*``median" + -0.208*``population" + 0.208*``copyright" + -0.200*``income" + 0.190*``fair" + -0.171*``census" + -0.168*``km" + -0.165*``households" \\
 \hline
 	5. & 117.243 & 0.317*``image" + -0.196*``players" + 0.190*``copyright" + 0.176*``median" + 0.174*``age" + 0.173*``fair" + 0.155*``income" + 0.144*``population" + -0.134*``football" + 0.129*``households" \\
 \hline
 	6. & 100.451 &  -0.504*``players" + -0.319*``football" + -0.284*``league" + -0.194*``footballers" + -0.141*``image" + -0.132*``season" + -0.117*``cup" + -0.113*``club" + -0.110*``baseball" + -0.103*``f" \\
 \hline
 	7. & 92.376 & 0.411*``album" + 0.275*``albums" + 0.217*``band" + 0.215*``song" + 0.184*``chart" + 0.164*``songs" + 0.160*``singles" + 0.149*``vocals" + 0.139*``guitar" + 0.129*``track" \\
 \hline
 	8. & 84.024 &  0.246*``wikipedia" + 0.183*``keep" + -0.179*``delete" + 0.167*``articles" + 0.153*``your" + 0.150*``my" + -0.141*``film" + 0.129*``we" + 0.123*``think" + 0.121*``user" \\
 \hline
 	9. & 79.548 &  \emph{word ``category'' in ten different languages (and their exotic un-\TeX-able scripts)} \\
 \hline
 	10. & 79.074 &  -0.587*``film" + -0.459*``films" + 0.129*``album" + 0.127*``station" + -0.121*``television" + -0.119*``poster" + -0.112*``directed" + -0.109*``actors" + 0.095*``railway" + -0.085*``movie" \\
 \hline
\end{tabular}

\end{document}